\documentclass[11pt]{amsart}

\newif\ifpdf
    \ifx\pdfoutput\undefined
    \pdffalse % we are not running PDFLaTeX
    \else
    \pdfoutput=1 % we are running PDFLaTeX
    \pdftrue
    \fi

    \ifpdf
    \usepackage[pdftex]{graphicx}
    \else
    \usepackage{graphicx}
    \fi
\usepackage{amsmath,amssymb,amsthm,epsfig}

\newcommand\A{{\mathbb H}^2}
\newcommand\R{\mathbb{R}}
\newcommand\Q{\mathbb{Q}}
\newcommand\N{\mathbb{N}}
\newcommand\Z{\mathbb{Z}}
\newcommand\HT{\A \times T_p}

\newcommand\F{\Z [\frac{1}{p}]}
\newcommand\PZp{PSL_2( \F)}

\newcommand\si{\sigma_{\infty}}

\newcommand\Qp{\mathbb{Q}_p}
\newcommand\PQp{PSL_2(\Qp)}
\newcommand\G{PSL_2(\R) \times PSL_2( \Qp)}
\newcommand\x{\alpha_0}

\newtheorem{theorem}{Theorem}[section]
\newtheorem{proposition}[theorem]{Proposition}
\newtheorem{lemma}[theorem]{Lemma}
\newtheorem{corollary}[theorem]{Corollary}

\begin{document}

\title{The Dehn function of $\PZp$}
\author{Jennifer Taback}
%\email{jtaback@math.albany.edu}
\thanks{The author would like to thank Jos\'e Burillo,
Benson Farb, David Fisher, Susan Hermiller and Sarah Rees for
useful mathematical conversations during the writing of this
paper, Kevin Wortman for corrections, and the referee for helpful
comments and great patience.}
%\institute{Department of Mathematics and Statistics\\ University at Albany \\
%  Albany, NY 12222}

\begin{abstract}
We show that $\PZp$ admits a combing with bounded asynchronous
width, and use this combing to show that $\PZp$ has an
exponential Dehn function.  As a corollary, $\PZp$
has solvable word problem and is not an automatic
 group.
\end{abstract}

\maketitle

\section{Introduction}
\label{sec:intro}

Let $<S|R>$ be a finite presentation of a group $G$, with $|S| <
\infty$. If $w \in F(S)$ is a word in the free group generated by
the set $S$ and $w=_G1$, then $w$ can be expressed as
$$w = \Pi_{i=1}^{A(w)} s_i R_i s_i^{-1}$$
for $s_i \in F(S)$ and $R_i \in R$, where $A(w)$ denotes the
minimal number of conjugates of relators necessary to express $w$.
The {\em Dehn function} $\delta: \N \rightarrow [0, \infty)$
 of the presentation $<S|R>$ is given by
$$\delta(n) = \max_{|w| \leq n} A(w)$$
where $|w|$ is the number of letters in $w$.

Define the following equivalence relation on the set of functions
\linebreak $f: \N \rightarrow [0, \infty)$.  If $f$ and $g$ are
two such functions, and there are constants $A,B,C,D,E \in \N$ so
that
$$f(n) \leq Ag(Bn+C) + Dn + E$$
for all $n \in \N$, then we write $f \prec g$.
If $f \prec g$ and $g \prec f$ then $f$ is equivalent to $g$, written
$f \sim g$.

If $f$ and $g$ are Dehn functions arising from presentations of
quasi-isometric groups, then $f \sim g$. Since two different
presentations for a group $G$ yield quasi-isometric Cayley graphs,
define the Dehn function of a finitely presented group $G$ to be
the equivalence class of the Dehn function arising from any of its
presentations. Thus the Dehn function of a finitely generated
group is independent of the choice of finite generating set.

The following fact is a consequence of the above equivalence
relation.   If a group $G$, in a given presentation, has a Dehn
function $f(n)$ which is, for example, exponential in $n$ (or
linear or polynomial), then a Dehn function arising from any other
presentation of $G$  will also be exponential (or linear or
polynomial, respectively) in $n$.

If $h: \N \rightarrow [0,  \infty)$ is a function satisfying
$max_{|w| \leq n} A(w) \leq h(n)$ for all $n$, then $h$ is called
an {\em isoperimetric function} for $G$ in the presentation
$\langle S | R \rangle$.  We again group these isoperimetric
functions into equivalence classes via the equivalence relation
defined above to obtain isoperimetric functions for $G$.  The Dehn
function is the minimal isoperimetric function for $G$.

One can also define the (geometric) Dehn function of a Riemannian manifold
$M$ as follows.
Let $f: S^1 \rightarrow M$ be a Lipschitz map and let $\gamma = f(S^1)$
with $l(\gamma) = length(\gamma)$.
By the Whitney Extension Theorem, the continuous extension of $f$
to $\hat f: D^2 \rightarrow M$ can also be chosen to be Lipschitz, with
the same constant.
We define the area $A(\gamma)$ of $\gamma$ to be the minimal area of all
Lipschitz discs bounded by $\gamma$.
The Dehn function of $M$ is defined by
$$\delta_M(x) = \sup_{l(\gamma) \leq x} A(\gamma).$$

In this paper, we prove that the group $\PZp$ has an exponential
Dehn function.  We do this using a geometric model $\Omega_p$ for
the group, which is described below. The Dehn function of
$\Omega_p$ is well defined since we can define Lipschitz loops in
$\Omega_p$ and compute their area.  It follows from \cite{BT} that
the Dehn functions of the group $\PZp$ and the space $\Omega_p$
are equivalent in the sense defined above.

The space $\Omega_p$ has an intricate  boundary consisting of
quasi-isometrically embedded copies of the solvable
Baumslag-Solitar group $BS(1,p^2)$, an example of the more general
solvable Baumslag-Solitar group $BS(1,n)$, for $n > 1$, which has
presentation
$$BS(1,n) = <a,b | aba^{-1} = b^n>.$$
It is well known that the Dehn function of $BS(1,n)$ is
exponential for $n > 1$. (\cite{G}, \cite{E+}) Using the boundary
components of $\Omega_p$, we easily obtain an exponential lower
bound for the Dehn function of $\Omega_p$, and thus of $\PZp$. To
determine that this Dehn function is actually exponential, we
construct a combing of $\PZp$ with bounded asynchronous width, and
apply a result of Bridson \cite{B} to obtain an exponential upper
bound.

\section{The geometry of $\PZp$}
\label{sec:geom}

\subsection{The geometric model}
\label{sec:omegap}

To construct a geometric model for $\PZp$, we will consider
$$PSL_2(\Z[ {\hbox{$\frac{1}{p}$}}]) \subset \G
$$
as the image of the diagonal map
$\eta$ given by $\eta(M) = (M,M)$ for $M \in \PZp$.
Viewed in this way, $\PZp$ is a lattice in the group $\G$, for any
prime $p$.

The group $\G$ acts by isometries on the space $\HT$, where $T_p$ is
the Bruhat-Tits tree associated to $PGL_2(\Qp)$ (see \cite{S} for the
construction of $T_p$).
  This follows from the fact that $PSL_2(\R)$ acts by isometries on $\A$ and $PSL_2(\Qp)$ acts by
isometries on $T_p$.
The tree $T_p$ is the regular $(p+1)$-valent tree, oriented so that
each vertex has one incoming and $p$ outgoing edges.

We will always view $\A$ as two dimensional hyperbolic space in the upper half space model.
Namely, $\A = \{ (x,y) | x \in \R, \ y > 0 \}$ with the metric $\frac{dx^2 + dy^2}{y^2}$.

The goal is to find a space $\Omega_p \subset \HT$ on which
$\PZp$ acts properly discontinuously and cocompactly by isometries.
The Milnor-Svarc criterion \cite{M} states that if a finitely generated
group $\Gamma$ acts
properly discontinuously and cocompactly by isometries on a space $X$,
then $\Gamma$ is quasi-isometric to $X$.
It follows that $\PZp$ and $\Omega_p$ are quasi-isometric.

The restriction of the action of $\G$ on $\HT$ to $\PZp$ is not
cocompact, although it is still by isometries and properly
discontinuous. Namely, the fundamental domain $D$ for the action
of $\PZp$ on $\HT$, when intersected with $\A$, is the same as the
fundamental domain for the action of $PSL_2(\Z)$ on $\A$, which is
unbounded in one direction. We now describe a method for
constructing a subspace of $\HT$ on which this action is
cocompact. (For a more complete construction of this subspace, we
refer the reader to \S $3$ of \cite{T}.)

Let $w$ be the intersection of the horocircle $y=h_0$ based at $\infty$
in $\A$ with the fundamental domain $D$.
Let $w'$ be a lift of $w$ to a horocircle $y=B$ in a fixed copy of $\A \subset \HT$.
The orbit of the segment $w'$ under $\PZp$ is a disjoint
collection ${\mathcal H}$ of horocircles based at $\Q \cup \{ \infty \} \subset
\partial_{\infty}\A$ in each copy of $\A \subset \HT$.

A {\em horosphere} $\sigma_{\alpha}$ of $\HT$ based at $\alpha \in
\Q \cup \{ \infty \}$ is defined to be a particular subset of
${\mathcal H}$, namely all horocircles in ${\mathcal H}$ based at
$\alpha$. For each copy of $\A \subset \HT$ there is a unique
horocircle in ${\mathcal H}$ lying in $\A$ which is based at
$\alpha$.

In order to have a ``connected'' picture of a horosphere,
add an edge $e$ between adjacent vertices $v_1$ and $v_2$ of $T_p$.
Let $h_1$ and $h_2$ be the horocircles in ${\mathcal H}$ based at
$\alpha$ in $\A \times \{v_1\}$ and $\A \times \{v_2\}$ respectively.
Extend the horosphere linearly in $\A \times e$ from $h_1$ to $h_2$,
obtaining a connected picture of a horosphere.

The space $\Omega_p$, where $\PZp$ acts properly
discontinuously and cocompactly by isometries, is defined to be
$\HT$ with the interiors of all the horospheres removed.  The interior
of a horosphere is the union of the interiors of the component
horocircles.

Below we will consider the horosphere $\si$ of $\Omega_p$ based at
$\infty$.
By construction, for any $\A \subset \HT$, the intersection $\si \cap \A$ is a horocircle
$h \subset \A$ based at $\infty$.
A matrix computation shows that these intersections are the
horocircles $\{ y=p^{2r}B \}$ for $r \in \R$, where $r$ changes
continuously, increasing along lines in $T_p$ and $B$ is fixed above.
Topologically $\si$, and thus any horosphere, is $\R \times T_p$.

The following theorem states that the product metric on $\HT$,
restricted to $\Omega_p$, is Lipschitz equivalent to the word
metric on $\PZp$.  In \S \ref{sec:combing} we use this metric to
define an asynchronous combing of $\Omega_p$.

\begin{theorem} \cite{LMR}
\label{thm:metric}
If $G$ is a semisimple Lie group of rank at least $2$ and $\Gamma$ is
an irreducible lattice in $G$ then, $d_R$ restricted to $\Gamma$ is
Lipschitz equivalent to $d_W$, where $d_W$ is the word metric on
$\Gamma$ and $d_R$ is the left invariant Riemannian metric on
$\Gamma$.
\end{theorem}

By construction, $\PZp$ acts properly discontinuously and cocompactly
by isometries on $\Omega_p$.

\subsection{The geometry of the horospheres}
\label{sec:horospheres}

The horospheres $\{ \sigma_{\alpha} \}$ of $\Omega_p$ also have
interesting geometric structure. The solvable Baumslag-Solitar
group $BS(1,n) = <a,b | aba^{-1} = b^n>$ for $n> 1$ acts properly
discontinuously and cocompactly by isometries on a metric
$2$-complex $X_n$ defined explicitly in \cite{FM}. This complex
$X_n$ is topologically $T_n \times \R$, where $T_n$ is a regular
$(n+1)$-valent tree, directed so that each vertex has $1$ incoming
edge and $n$ outgoing edges.

\begin{figure}
\begin{center}
\includegraphics[width=1.5in]{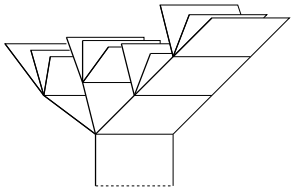}
\caption{The geometric model of the solvable Baumslag-Solitar group $BS(1,3)$.}
\label{fig:bs}
\end{center}
\end{figure}

We can describe $X_n$ metrically after fixing a height function
$h$ on the tree $T_p$.  A {\em height function} on $T_p$ is a
continuous function $h:T_p \rightarrow \R$ which maps each
oriented edge of $T_p$ homeomorphically onto an oriented interval
of a given length $d$. Fix a basepoint $t_0$ for $T_p$ with height
$0$. Then each point $t$ of $T_p$ is assigned a height $h(t)$; for
example, a vertex $t$ connected to $t_0$ by an edge leaving $t_0$
has height $d$.  If $l \subset T_p$ is a line on which $h$ is
strictly increasing, then $l \times \R \subset X_p$ is metrically
a hyperbolic plane.

The map sending the generators $a$ and $b$ of $BS(1,p^2)$ to the matrices
$A = \left( \begin{array}{cc} p & 0 \\ 0 & \frac{1}{p}
\end{array} \right)$ and
$B= \left( \begin{array}{cc} 1 & 1 \\ 0 & 1
\end{array} \right)$ respectively is a homomorphism of $BS(1,p^2)$ into the group $\PQp$.
The orbit of the segment $w' \subset \HT$ (defined in \S
\ref{sec:omegap}) under the group generated by the matrices $A$
and $B$ is exactly the horosphere $\si$ of $\Omega_p$. (\cite{T},
\S 3.5)

One checks that the horosphere $\si$ of $\Omega_p$ (and hence any
horosphere of $\Omega_p$) is a quasi-isometrically embedded copy
of the complex $X_{p^2}$ corresponding to $BS(1,p^2)$. (\cite{T},
\S 3.5) The Dehn function of $BS(1,n)$ is known to be exponential
for $n
> 1$. (See, e.g. \cite{E+} or \cite{G}.)

\section{A combing of $\Omega_p$}
\label{sec:combing}

Below we construct a combing of $\Omega_p$ which we show to have
bounded asynchronous width.  We will eventually use this combing
to obtain an exponential upper bound on the Dehn function of
$\PZp$. We follow \cite{B} in our definitions.

Let $G$ be a group with a finite generating set $X = \{ x_1
, \cdots x_k\}$, and let $X^*$ denote the set of all
words in the letters of $X$ (and their inverses).
A {\em combing} of $G$ is a choice, for
each $g \in G$, of a word in $X^*$ representing $G$, i.e. a section of the
natural map $X^* \rightarrow G$.  Bridson
does not require his combing paths to be quasigeodesics, although
later we will add some conditions on the length of the
combing paths in order to obtain the results that we need.

One can view the combing geometrically as a set of paths in the Cayley
graph $\Gamma(G,X)$: for each $g \in G$, the combing gives a unique choice of
path from the identity to $g$.
We keep track of the {\em width} between these paths, which is a bound on
the distance between points moving at unit speed along combing lines
which end at points distance 1 apart, as well as the
length of the paths.  If $\sigma: G \rightarrow X^*$ is a combing,
let $\sigma_g(t)$ denote the (parametrized) combing path from the identity
to $g \in G$.  The length function for the combing is defined to be
$$L(n) = \max \{ length(\sigma_g) | d(id, g) \leq n \}.$$

For this paper, we are interested in {\em asynchronous combings}, in which
we are allowed to consider the combing paths traversed at different speeds.
Rather than observing the width of the combing, we instead consider the
asynchronous width, which is defined as follows.

Fix a combing of the finitely generated group $G$, namely $\sigma:
G \rightarrow X^*$. We first consider the set of all possible
reparametrizations of our combing paths, i.e. the set of functions
$$
R = \{ \rho: \N \rightarrow \N | \rho(0) = 0, \ \rho(n+1) \in
 [\rho(n),\rho(n)+1] \ \forall n, \ \rho \text{ unbounded} \}.$$
Given two elements $g,h \in G$, we choose functions $\rho$ and
$\rho'$ in $R$ with which we reparametrize the combing paths
$\sigma_g$ and $\sigma_h$, respectively, and we consider the width
between these paths. We then minimize this width over all possible
reparametrization functions in the set $R$.  In summary, we define
$$ D_{\sigma}(g,h) = min_{\rho, \ \rho' \in R}
\{ \max_{t \in \N} \{ d(\sigma_g(\rho(t)),\sigma_h(\rho'(t)))\} \}.$$

Then the {\em asynchronous width} of $\sigma$ is defined to be the function
$$
\phi(n) = \max \{D_{\sigma}(g,h) | d(id,g) \leq n, \ d(id, h) \leq n,
\ d(g,h) = 1 \}.$$
All distances above are computed using the word metric in the Cayley graph
$\Gamma(G,X)$.

Let $\Omega$ be a metric space with a fixed basepoint $x_0$. We
can define a combing of $\Omega$ as a choice of path from $x_0$ to
each point $x \in \Omega$. We define the length function $L(n)$
and asynchronous width function analogously, replacing any
instance of the identity of the group with the basepoint $x_0$.
The reparametrization functions in the set $R$ must now be
continuous functions from $\R^* = \R^+ \cup \{ 0 \} $ to $\R^*$,
satisfying the same conditions as the discrete functions in $R$
defined above. In particular, we now construct a combing of the
space $\Omega_p$ defined in \S \ref{sec:geom}.

We will show that the combing constructed below of
$\Omega_p$ has bounded asynchronous width.  Since $\PZp$ and
$\Omega_p$ are quasi-isometric, we explain below that the
associated combing of $\PZp$ will also have bounded
asynchronous width.

It will also follow from the construction of the combing paths in
$\PZp$ that the length of $\sigma_g$ is at most exponential in
$d(id, g)$, for $g \in \PZp$.  This fact, together with the
bounded asynchronous width of the combing, allows us to apply a
theorem of Bridson which guarantees an exponential upper bound for
the Dehn function of $\PZp$.

\subsection{Constructing the combing paths}
\label{sec:paths}

Fix a height function $h$ on $T_p$, and a basepoint $t_0 \in T_p$
at height $0$.  We use the coordinates $(x,y,t)$ on $\HT$, where
$(x,y)$, with $y
> 0$, denotes a point in $\A$ and $t \in T_p$. Let $\alpha =
(x',y',t')$ be any point in $\Omega_p$. We now describe a
canonical path from a fixed basepoint $\x = (x_0,y_0,t_0) \in
\Omega_p$ to $\alpha$. Consider the following pair of paths in
$\HT$.  If $\xi_1$ is the unique path in $T_p$ lying between $t_0$
and $t'$ and $\xi_2$ is the unique hyperbolic geodesic between the
points $(x_0,y_0)$ and $(x',y')$, then let
$$\gamma_1= \{ (x_0,y_0,t) | t \in \xi_1 \} \text{ and } \gamma_2= \{
(x,y,t') | (x,y) \in \xi_2 \}.$$ So $\gamma_1$ lies entirely in
the tree factor and $\gamma_2$ lies entirely in an $\A$ factor.
The composition $\gamma = \gamma_2 \circ \gamma_1$ is a path from
$\x$ to $\alpha$, in which the path $\gamma_1$ is followed by
$\gamma_2$. Following theorem \ref{thm:metric}, we use the product
metric on $\Omega_p$ as a subspace of $\HT$. Thus $d(\alpha_0,
\alpha) = length(\gamma_1) + length (\gamma_2)$. The path $\gamma$
may not lie completely in $\Omega_p$.  We adapt it below so that
it will be contained in $\Omega_p$.

\subsubsection{Adapting the path in the tree factor}
\label{sec:tree}

Suppose that $I$ is a maximal connected closed interval of $\gamma_1$
of length $k$
not lying in $\Omega_p$, but inside a horosphere $\sigma$.  Without loss of
generality we assume that $\sigma = \si$.

Let $v_1=(x_0,y_0,t_1)$ and $v_2=(x_0,y_0,t_2)$ be the endpoints
of $I$. Either both $v_i$ lie in $\si$ or we have $v_1 \in \si$
and $v_2$ in the interior of $\si$. In the first case, when both
endpoints are in $\si$, let $\eta_1$ be the unique path in $T_p$
between $t_1$ and $t_2$. We replace $I$ by the path $I'$ in $\si$
lying between $v_1$ and $v_2$. Recall from \S \ref{sec:omegap}
above and \S 3 of \cite{T} that points in $\si$ can be described
by coordinates $(x,p^{2h(t)}B,t)$ where $h$ is a fixed height
function on $T_p$, $x \in \R$, with $t \in T_p$ and $B> 1$ is a
fixed constant. Thus the path $I'$ has the form
$(x_0,p^{2h(t)}B,t)$ for $t \in \eta_1$.

If $v_1 \in \si$ and $v_2$ lies in the interior of $\si$, then let
$\eta_1$ be as above, and replace $I$ by the path $I' = \{
(x_0,p^{2h(t)}B,t) | t \in \eta_1 \}$. Notice that in this case,
$I'$ has one distinct endpoint from $I$, namely
$(x_0,p^{2h(t')}B,t')$, which is different from $v_2$. Since $I$
was chosen to be a maximal connected interval of $\gamma_1$ not
contained in $\Omega_p$, the only way this case occurs is if $v_2$
is the endpoint of $\gamma_1$, i.e. $v_2 = (x_0,y_0,t')$. Let
$v_2'$ denote this new endpoint. When we adapt $\gamma_2$ below we
will ensure that the path from $\alpha_0$ to $\alpha$ is
connected.

In either case the path $I'$ is composed of segments of vertical geodesics in the
hyperbolic planes of $\si$.
If the length of $I$ was $l$, the length of $I'$ is at most $K'l$
where $K'$ is the fixed quasi-isometry constant obtained when
considering the quasi-isometric embedding of the horospheres into
$\Omega_p$.  \cite{T}

For each such interval $I \subset \gamma_1$, replace $I$ by $I'$ constructed
as above, obtaining a new path $\gamma_1'$, satisfying
$length(\gamma_1') \leq K' length(\gamma_1)$.
Note that $\gamma_1'$ may have a different (non-basepoint) endpoint
$v_2'$ than $\gamma_1$.

\subsubsection{Adapting the path in the hyperbolic factor}
\label{sec:hyperbolic}
Now suppose that $J$ is a maximal connected closed interval of $\gamma_2$
of length $k$
not lying in $\Omega_p$, but inside a horosphere $\sigma$.
If there is no such interval $J$ then $\sigma_{\alpha} = \gamma_2 \circ \gamma_1'$
is the desired path from $\alpha_0$ to $\alpha$, completely contained in
$\Omega_p$.

Now assume that such an interval $J$ exists.
Without loss of
generality we assume that $\sigma = \si$ and that the endpoints of $J$
are $v_3$ and $v_4$.

Suppose that both $v_3$ and $v_4$ lie on $\si$.  In this case we
replace $J$ with the horocyclic segment $J'$ lying in $\si$
between the endpoints.
It is well known that downward projection along vertical geodesics
in $\A$ onto a horosphere increases length exponentially. Thus the
length of $J$ is increased exponentially.

If this is not the case, then $v_3=v_2$ lies in the interior of $\si$
and is the common endpoint of
$\gamma_1$ and $\gamma_2$.
Recall that $v_2 = (x_0,y_0,t')$ and define $v_2' = (x_0,p^{2h(t')}B,t')$ to be the projection of $v_2$ onto $\si$ along vertical geodesics.
Replace $J$ with the the horocyclic segment $J'$ of $\si$ between $v_2'$
and $v_4=(x_4,p^{2h(t')}B,t')$, namely the path
$$J' = \{ (x,p^{2h(t')}B,t') | x \in [x_0,x_4]\}.$$

Replacing each interval $J$ of this form by the corresponding
interval $J'$ constructed as above, we obtain a new path
$\gamma_2'$, whose length satisfies $length(\gamma_2') \leq
e^{length(\gamma_2)}$. By construction, $\gamma_1'$ and
$\gamma_2'$ have a common endpoint, so that $\gamma_2 \circ
\gamma_1$ is a path from $\x$ to $\alpha$, whose length is at most
exponential in $length(\gamma)$.  We define the combing path
$\sigma_{\alpha}$ to be $\gamma_2 \circ \gamma_1$.

We summarize conditions on the length of the combing paths in the
following lemma.

\begin{lemma}
\label{lemma:length} Let $\sigma$ denote the combing of $\Omega_p$
defined above. Then the length function of $\sigma$ satisfies an
exponential upper bound, i.e. $L(n) \leq e^{n}$.
\end{lemma}

\subsubsection{Bounded asynchronous width}
We now prove that the asynchronous width of the combing of
$\Omega_p$ defined above is bounded.  Let $\alpha_1$ and
$\alpha_2$ denote points in $\Omega_p$ at distance $1$ from each
other, and let $\sigma_{\alpha_i}$ denote the combing path from a
fixed basepoint $\alpha_0$ of $\Omega_p$ to $\alpha_i$ for
$i=1,2$. We keep the notation $\gamma_2 \circ \gamma_1$ for the
combing paths, so that we can tell which part of each path lies in
the hyperbolic factor ($\gamma_2$), and which in the tree factor
($\gamma_1$).

We begin with two lemmas which prove that geodesics in $\A$ emanating from a common point and ending one unit apart have bounded asynchronous width.

We define some notation used in the lemmas below. Let $(x_1,y_1)$
and $(x_2,y_2)$ be points in $\A$ at distance $1$ from each other.
By applying a hyperbolic isometry $\phi$, we may assume that $(x_1,y_1)$
and $(x_2,y_2)$ lie on a vertical geodesic $x=n$. Let $z=(z_1,z_2)
\in \A$ be a point distinct from $(x_1,y_1)$ and $(x_2,y_2)$, not
lying on the same vertical geodesic, and let $\chi_i$ be the
geodesic segment connecting $z$ and $(x_i,y_i)$ for $i = 1,2$. It
will be useful in the arguments below to express $\chi_1$ and
$\chi_2$ in the coordinates $(x,y)$ of $\R^2$.  Namely, after applying the hyperbolic isometry $\phi$ to ensure that $(x_1,y_1)$ and $(x_2,y_2)$ lie on a common geodesic, we will assume that
$\chi_1$ lies on the circle with equation $(x-a)^2 + y^2 = R_1^2$
and $\chi_2$ on the circle with equation $(x+a)^2 + y^2 = R_2^2$,
for some $a > 0$.  Also assume that each $\chi_i$ is parameterized by $s \in [0,
length(\chi_i)]$, and that $length(\chi_2) > length(\chi_1)$.

\begin{lemma}
\label{lemma:reparam} Let $(x_i,y_i)$ and $\chi_i$ be as above,
for $i=1,2$. Then there exists a function $\rho' \in R$ which is
strictly increasing with the property that $\chi_2(s)$ and
$\chi_1(\rho'(s))$ lie on a vertical geodesic in $\A$, for $s \in [0, length(\chi_2)]$.
\end{lemma}

\begin{proof}
It is clear that there exists such a function $\rho'$ so that the
points in question lie on the same vertical geodesic; we now show
that this function is increasing. We first note that as functions
of $s$, the parameterizations $\chi_1$ and $\chi_2$ are strictly
increasing, because the parameter $s$ represents path length.  We
will write $\chi_i(s_1) < \chi_i(s_2)$ to mean that the point
$\chi_i(s_2)$ is further along the path $\chi_i$ (beginning at
$z$) than $\chi_i(s_1)$.

Suppose $\rho'(x)$ was not increasing, so there are points $s_1, \
s_2 \in [0,length(\chi_1)]$ with $s_1 < s_2$ but $\rho'(s_1) \geq
\rho'(s_2)$.  Then $\chi_1(\rho'(s_1)) \geq \chi_1(\rho'(s_2))$.
However, $\chi_2(s_1) < \chi_2(s_2)$, and $\chi_2(s_i)$ and
$\chi_1(\rho'(s_i))$ lie on a common vertical geodesic, a
contradiction.
\end{proof}

\begin{lemma}
\label{lemma:Kfellowtraveller} Let $(x_i,y_i)$ and $\chi_i$ for
$i=1,2$ be as above, and $\rho'(x) \in R$ the function guaranteed
in lemma \ref{lemma:reparam}.  Then for any $s \in [0,
length(\chi_2)]$, we have $d_{\A} (\chi_1(\rho'(s)), \chi_2(s))
\leq 1$, where $d_{\A}$ represents distance in the hyperbolic
plane $\A$.
\end{lemma}

\begin{proof}
Using the equations given above for the circles containing
$\chi_1$ and $\chi_2$, we express the hyperbolic distance between
the points $\chi_1(\rho'(s))$ and $\chi_2(s)$ which lie on a
vertical geodesic as a function $f(s)$ of $s$.  Namely,
$$f(s) = \ln \left( \frac{\sqrt{R_1^2 - (t-a)^2}}{\sqrt{R_2^2 - (t+a)^2}}
\right) = \frac{1}{2} \ln \left( \frac{R_1^2 - (t-a)^2}{R_2^2 -
(t+a)^2} \right).$$ To prove the lemma we show that this function
is strictly increasing. Then, since $$f(l(\chi_2)) =
d_{\A}(\chi_1(\rho(length(\chi_1)), \chi_2(length(\chi_2))) =
d_{\A} (x_1,x_2) \leq 1,$$ the inequality in the lemma follows.

It is sufficient to show that the function $$g(x) = \ln \left(
\frac{R_1^2 - (x-a)^2}{R_2^2 - (x+a)^2} \right) = \ln ( R_1^2 -
(x-a)^2) - \ln (R_2^2 - (x+a)^2)$$ is increasing.  Computing
$g'(s)$ we see that
$$g'(s) = \frac{-2(x-a)}{R_1^2 - (x-a)^2} - \frac{-2(x+a)}{R_2^2 - (x+a)^2}
\geq \frac{4a} {R_1^2 - (x-a)^2} >0,$$ since $a >0$.  Thus $g(s)$,
and more importantly $f(s)$, are strictly increasing, as desired.
\end{proof}

In the following corollary, we no longer assume that the points $(x_1,y_1)$ and $(x_2,y_2)$ lie on a vertical geodesic.

\begin{corollary}
\label{cor:rho}
Let $(x_i,y_i) \in \A$ for $i=1,2$ with $d_{\A}((x_1,y_1),(x_2,y_2)) \leq 1$, and $z \in \A$.
Let $\chi_i$ denote the hyperbolic geodesic from $z$ to $(x_i,y_i)$.  Then there is a
reparametrization function $\rho \in R$ so that $d_{\A}(\chi_1(\rho(s)), \chi_2(s))
\leq 1$ for all $s \in [0, length(\chi_2)]$.
\end{corollary}

\begin{proof}
Let $\phi$ denote a hyperbolic isometry such that $\phi(x_1,y_1)$ and $\phi(x_2,y_2)$ lie on a vertical geodesic.
Let $\rho'$ be the reparametrization function guaranteed in lemma \ref{lemma:reparam}.  Then $\rho=\phi^{-1} \circ \rho' \circ \phi$ satisfies the conditions of the corollary.
\end{proof}

Returning to $\Omega_p$, let $\alpha_1, \  \alpha_2 \in \Omega_p$
with $d(\alpha_1,\alpha_2) \leq 1$. Let $\sigma_{\alpha_i}$ for
$i=1,2$ denote the combing path contained in $\Omega_p$ described
in \S \ref{sec:paths} from the basepoint $\x$ to $\alpha_i$. Then
$\sigma_{\alpha_1} = \gamma_2 \circ \gamma_1$ and
$\sigma_{\alpha_2} = \xi_2 \circ \xi_1$ where $\tilde \gamma_2$
and $\tilde \xi_2$ are each contained entirely in the hyperbolic
factor and $\tilde \gamma_1$ and $\tilde \xi_1$ are contained
entirely in the tree factor before they are adapted to lie in
$\Omega_p$.  We keep the notation $\tilde \gamma_i$ and $\tilde
\xi_i$ to denote the combing paths before they have been adapted
to lie in $\Omega_p$.

Let $\pi$ denote downward projection along vertical geodesics onto $\si$ in any $ \A \subset \HT$.  This is the projection used to adapt the initial hyperbolic geodesics into combing paths as described above.

Our goal is to prove the following proposition.

\begin{proposition}
\label{prop:boundedwidth}
The combing of $\Omega_p$ defined above
has bounded asynchronous width.
\end{proposition}

It is clear that since $d(\alpha_1,\alpha_2) \leq 1$, the pieces
$\tilde \gamma_1$ and $\tilde \xi_1$ of the combing paths which
lie in the tree factor can be reparametrized by a function $\psi
\in R$ so that $d(\tilde \gamma_1(\psi(s)), \tilde \xi_1(s) )\leq
1$ for $s \in [0, length(\tilde \xi_1)]$.  It immediately follows
that $d(\gamma_1(\psi(s)), \xi_1(s)) \leq 2 \log p$.  We are
therefore concerned with showing the analogous condition for the
pieces of the combing paths which lie in the hyperbolic factor. We
do this via a reparametrization of the unadapted combing paths
$\tilde \gamma_2$ and $\tilde \xi_2 $  which is inherited by the
adapted combing paths.  We divide the proof into several cases. We
first consider the case when at least one of $\tilde \gamma_2$ and
$\tilde \xi_2$ is completely contained in $\Omega_p$.  If both
paths are completely contained in $\Omega_p$, then $\tilde
\gamma_2 = \gamma_2$, and $\tilde \xi_2 = \xi_2$, and the bounded
asynchronous width is immediate from corollary \ref{cor:rho}.

To obtain bounded asynchronous width when $\tilde \gamma_2$ intersects $\si$, we need only the following basic hyperbolic geometry lemma.

\begin{lemma}
\label{lemma:hyperbolic} Let $m=(x_1,y_1)$ and $n=(x_2,y_2)$ be
points in $\A$ so that $d_{\A}(m,n) \leq 1$, and let $y=h$ be any
horocircle intersecting the geodesic segment $[m,n]$ at a point $k
\neq m,n$ labelled so that $y_1 < h < y_2$. Let $n'$ be the
projection $\pi(n)$ of $n$ onto $y=h$ along vertical geodesics.
Then $d_{\A}(m,n') \leq 2$.
\end{lemma}

\begin{proof}
The lemma follows from the fact that the shortest distance between a point and a horocircle in $\A$ is measured along a vertical geodesic, and the triangle inequality.
\end{proof}

The following corollary is immediate.

\begin{corollary}
Suppose that $\gamma_2$ and $\xi_2$ are components of combing paths $\sigma_{\alpha_1}$ and $\sigma_{\alpha_2}$, respectively, with $d(\alpha_1,\alpha_2) \leq 1$, as defined above, arising from the hyperbolic geodesics $\tilde \gamma_2$ and $\tilde \xi_2$, where $\tilde \xi_2 \subset \Omega_p$ but $\tilde \gamma_2 \cap \si \neq \emptyset$.
Then for a suitable reparametrization function $\rho \in R$, we have $d_{\A}(\xi_2(t),\gamma_2(\rho(t)))$ uniformly bounded, for $t \in [0, length(\xi_2)]$.
\end{corollary}

\begin{figure}
\begin{center}
\includegraphics[width=4in]{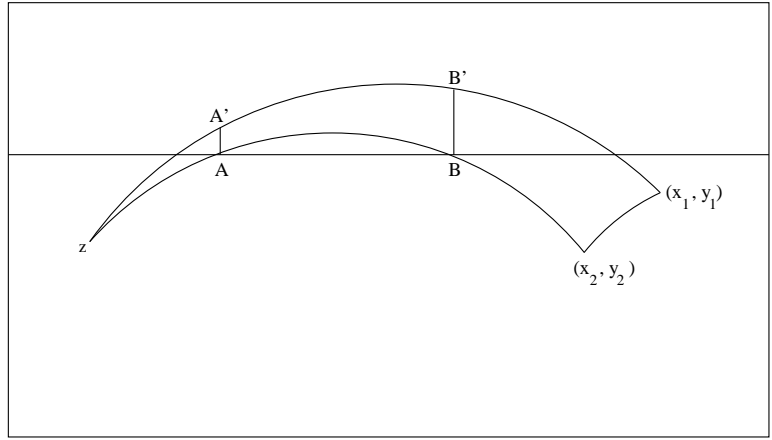}
\caption{A hyperbolic plane $\A \subset \HT$ with the horocircle
$\si \cap \A$ depicted as well.  The hyperbolic geodesic
containing the points $A$ and $B$ is $\tilde \xi_2$ in the
arguments below, from $z$ to $\alpha_2 = (x_2,y_2)$.  The other
hyperbolic geodesic depicted above is $\tilde \gamma_2$ from $z$
to $\alpha_1 = (x_1,y_1)$. \label{fig:intersection}}
\end{center}
\end{figure}
\begin{proof}
Reparametrize $\tilde \gamma_2$ using the function $\rho$ produced
above in corollary \ref{cor:rho}, and project $\tilde \gamma_2$
via $\pi$ to obtain $\gamma_2$. Lemma \ref{lemma:hyperbolic} then
asserts that $d_{\A}(\gamma_1(t),\gamma_2(\rho(t)) )\leq 2$.
\end{proof}

We now consider the case when both $\tilde \gamma_2$ and $\tilde
\xi_2$ have nontrivial intersection with $\si$, as depicted in figure \ref{fig:intersection}.  In this case we
define completely new parametrizations of both paths in order to
show that the asynchronous width is bounded.  We do this below in
several cases.  Many of the cases involve the original
parametrization $\rho$ of corollary \ref{cor:rho}.

We define a parametrization $\omega_1$ of $\tilde \xi_2$  and
$\omega_2$ of $\tilde \gamma_2$.   Note that we have not yet
projected the paths so that they lie in the space $\Omega_p$, and
thus $\tilde \xi_2$ and $\tilde \gamma_2$ are hyperbolic
geodesics.   We begin with a few distinguished points.  Let $A$
and $B$ be the points of intersection of $\tilde \xi_2$ with
$\si$, where $A$ is closer to $z$ than $B$.  Using the initial parametrization of $\tilde \xi_2$ by $t$ and of
$\tilde \gamma_2$ by $\rho(t)$, for $t \in [0,length(\tilde \xi_2)]$,  let $t_A$ denote the parameter
value so that $\tilde \xi_2(t_A) = A$, and $t_B$ the parameter
value so that $\tilde \xi_2(t_B) = B$.

Now let $\bar{A}  = \tilde \gamma_2(\rho(t_A))$, and $\bar{B} = \tilde \gamma_2(\rho(t_B))$.  Let $A'$
(resp. $B'$) be the point of $\tilde \gamma_2$ lying on the vertical geodesic containing $A$ (resp. $B$).
Define $t_{A'} \in [0, length(\tilde \xi_2)]$ to  be the parameter value so that $\tilde \gamma_2(\rho(t_{A'})) = A'$.  Define $t_{B'}$ analogously.

We now define parametrizations $\omega_1$ of $\tilde
\xi_2$ and $\omega_2$ of $\tilde \gamma_2$ which will be used to
show that the adapted paths $\xi_2$ and $\gamma_2$, which are contained in $\Omega_p$, have bounded
asynchronous width.  We again use the notation $M<N$ if $M$ and
$N$ are two points on either $\tilde \xi_2$ or $\tilde \gamma_2$,
and $M$ is closer in hyperbolic distance to $z$ than $N$.

We will make use of the following function $v(t)$ in the
parametrizations defined below. Suppose that for some interval
$I \subset [0, length(\tilde \xi_2)]$, we parametrize $\tilde
\xi_2$ by $\tilde \xi_2(t)$ for $t \in I$. Let $v(t)$ for $t \in
I$ be defined so that $\tilde \gamma_2(v(t))$ is the unique point
of $\tilde \gamma_2$ which lies on the same vertical geodesic in
$\A$ as $\tilde \xi_2(t)$.

The idea behind the parametrizations is the following.  Keeping in mind that to form the ultimate combing paths we must project $\tilde \xi_2$ and $\tilde \gamma_2$ via the projection $\pi$, define $\omega_1$ and $\omega_2$ so that between the points $A$ and $B$, the same parameter value corresponds to points on $\tilde \xi_2$ and $\tilde \gamma_2$ lying on a common vertical geodesic.  The images of these points under $\pi$ will be identical.   The remainder of the definitions of $\omega_1$ and $\omega_2$ ensures that after the projection $\pi$, the bounded asynchronous width condition will be satisfied for the rest of the paths.   The parametrizations $\omega_1$ and $\omega_2$ are defined as follows.

We first define the parametrization of both paths from the point z to the points $B$ and $B'$, respectively.
\begin{enumerate}
\item
If $\bar{A} > A'$, then define $\omega_1(t) = t$ for $t \in [0,t_B]$, and $\omega_2(t) = v(t)$ for $t \in [0,t_B]$.
\item
If $\bar{A} \leq A'$ begin by defining
$\omega_1(t) = \left\{ \begin{array}{ll} t & t \in [0,t_A] \\
t_A  & t \in [t_A,t_{A'}] \end{array}\right.$ and $\omega_2(t) =
\rho(t)$ for $t \in [0,t_{A'}]$.
We then use $\omega_1$ to rescale the interval $[t_{A'},t_B]$ so that at the endpoints of the interval, we have $\omega_1(t_{A'}) = t_A$ and $\omega_1(t_B) = t_B$.  For $s \in [t_{A'},t_B]$, define $\omega_2(s) = v(\omega_1(s))$.

\medskip
\noindent
We now continue the parametrizations to the end of each path.
\item
If $\bar{B} \leq B'$, then $\omega_1(t) = t$ for $t \in [t_B,length(\tilde \xi_2)]$, and
$\omega_2(t) = \left\{ \begin{array}{ll} t_{B'} & t \in [t_B,t_{B'}] \\
\rho(t)  & t \in [t_{B'},length(\tilde \xi_2)] \end{array}\right.$.

\item
If $\bar{B} > B'$, we must subdivide the interval $[t_B,
length(\tilde \xi_2)]$ in order to complete the parametrization
functions.  According to the definition of $\omega_2$, we have
that $\tilde \gamma_2(\omega_2(t_B)) = B'$.  We introduce an
intermediate point $d$ with $t_B < d < length(\tilde \xi_2)$, and
scale as follows.

We use a function $f$ to scale $[t_B,d]$ for the definition of
$\omega_2$, so that the following conditions hold. At the
endpoints of the interval, we have that $\omega_2(t_B) = f(t_B) =
\rho(t_{B'})$ and $\omega_2(d) = f(d) = \rho(t_{B})$. On the interior of the
interval, we use the definition $\omega_2(t) = \rho(f(t))$. We
define $\omega_1(t) = t_B$ for $t \in [t_B,d]$.

On the interval $[d, length(\tilde \xi_2)]$, we again use a
scaling function $g$, subject to the following conditions.  On the
endpoints, we require that $\omega_1(d) = g(d) = t_B$ and
$\omega_1(length(\tilde \xi_2)) = g( length(\tilde \xi_2))= length(\tilde \xi_2)$. Then
define $\omega_1(t) = g(t)$, and $\omega_2(t) = \rho(g(t))$ for $t \in [d,length(\tilde \xi_2)]$ .
\end{enumerate}

The parameterizations $\omega_1$ and $\omega_2$ apply to the {\em unadapted} paths $\tilde \xi_2$ and $\tilde \gamma_2$. We are interested in showing that the paths $\xi_2$ and $\gamma_2$ which lie in the space $\Omega_p$ have bounded asynchronous width.   These paths will inherit the parametrizations $\omega_1$ and $\omega_2$, respectively.

\begin{lemma}
\label{lemma:bounded}
Suppose that $\gamma_2$ and $\xi_2$ are components of combing paths $\sigma_{\alpha_1}$ and $\sigma_{\alpha_2}$, with $\alpha_1$ and $\alpha_2$ lying in a common hyperbolic plane, and $d(\alpha_1,\alpha_2) \leq 1$, as defined above, arising from the hyperbolic geodesics $\tilde \gamma_2$ and $\tilde \xi_2$, neither of which is completely contained in $\Omega_p$.
Then there are parametrization functions $\omega_1, \ \omega_2 \in R$ so that $d_{\A}(\xi_2(\omega_1(t)),\gamma_2(\omega_2(t))$ is uniformly bounded, for $t \in [0, length(\tilde \xi_2)]$.
\end{lemma}

\begin{proof}
We show that the parametrization functions $\omega_1$ and $\omega_2$ defined above satisfy the condition of the lemma.
We go through the proof according to the cases of the parameterizations above, showing in each case that the distance between the points is uniformly bounded.
\begin{enumerate}

\item Suppose that $\bar{A} > A'$.  Notice that $d_{\A}(A', \si)
\leq d_{\A}(\bar{A},\si) \leq d_{\A}(\bar{A},A) \leq 1$. Thus if we restrict to the initial
segments of $\tilde \xi_2$ and $\tilde \gamma_2$ from $z$ to $A$
and $z$ to $A'$ respectively, we see that we have a pair of
geodesics satisfying the conditions of lemma
\ref{lemma:Kfellowtraveller}. Thus we may conclude that points on
these paths lying on a common vertical geodesic are at most
distance $1$ apart as well.  It is then clear that the projection
$\pi$ along vertical geodesics can only decrease this distance.

For $t \in [t_A,t_B]$, we have that
$$\xi_2(\omega_1(t)) = \pi(\tilde \xi_2(\omega_1(t)))= \pi(\tilde \xi_2(t)) = \pi(\tilde \gamma_2(v(t))) = \pi(\tilde \gamma_2(\omega_2(t))) = \gamma_2(\omega_2(t))$$
and thus the distance between the points in question is $0$.

\item
Suppose that $\bar{A} \leq A'$.  By construction, $\tilde \xi_2([0,t_A]) = \xi_2([0,t_A]) \subset \Omega_p$.  Thus for any point $\gamma_2(\omega_2(t))$ for $t \in [0,t_A]$, either lemma  \ref{lemma:hyperbolic} or the definition of $\rho$, implies that $d_{\A}(\xi_2(\omega_1(t)),\gamma_2(\omega_2(t))) \leq 2$.

First consider the case when $\bar{A} \notin \Omega_p$.  If $t \in
[t_A,t_{A'}]$, it is clear that under the projection $\pi$, we
have $d_{\A} (A, \pi(\tilde \gamma_2(\omega_2(t)))) \leq d_{\A}(A,
\pi(\bar{A}))$. Then $\pi(\bar{A}), \ A, \ \text{and } \bar{A}$
form a geodesic triangle with $d_{\A}(\pi(\bar{A}),\bar{A}) \leq
d_{\A}(A, \bar{A}) \leq 1$ and thus $d_{\A}(\pi(\bar{A}),A) \leq
2$ forcing the distance between these points when measured along
the horocircle to be at most $e^2$.

If $\bar{A} \in \Omega_p$ and $\tilde \gamma_2(\omega_2(t)) \in \Omega_p$ for some $t \in [t_A,t_{A'}]$, the triangle inequality again shows that $d_{\A}(A, \gamma_2(\omega_2(t))) \leq 2$.  The triangle inequality also implies that if $D$ is the left most point of intersection of $\si$ with $\tilde \gamma_2$, then $d_{\A}(A,D) \leq e^2$.  The above argument shows that $d_{\A}(A, \gamma_2(\omega_2(t))) \leq e^2$ when $\tilde \gamma_2(\omega_2(t)) \notin \Omega_p$ and $t \in [t_A,t_{A'}]$.

For $t \in [t_{A'},t_B]$, we have that
$$\xi_2(\omega_1(t)) = \pi(\tilde \xi_2(\omega_1(t)))=\pi(\tilde \gamma_2(v(\omega_1(t)))) = \gamma_2(v(\omega_1(t))) = \gamma_2(\omega_2(t))$$
and thus the distance between the points in question is $0$.
\item
Suppose that $\bar{B} \leq B'$.  For $t \in [t_B,t_{B'}]$, it follows from the triangle inequality that $d_{\A}(\xi_2(t), \gamma_2(t_{B'}))  = d_{\A}(\xi_2(t),  B) \leq 2$.

For $t \in [t_{B'}, length(\tilde \xi_2)]$, we see that the geodesic $[\xi_2(t_{B'}) , \alpha_2]$ is contained in $\Omega_p$.  Lemma \ref{lemma:hyperbolic} and the definition of $\rho$ then imply that $d_{\A}(\xi_2(\omega_1(t)),\gamma_2(\omega_2(t))) \leq 2$ for $t$ in this interval.

\item
Suppose that $\bar{B} > B'$, and consider $t \in [t_B,d]$.  If $\bar{B} \in \si$, it follows from the above arguments that the distance $d_{\A}(\xi_2(t_B),\gamma_2(\omega_2(t))) =  d_{\A}(B,\gamma_2(\omega_2(t)))  \leq e^2$.  If $\bar{B} \in \Omega_p$, we must supply an additional argument. Let $D$ be the final point of intersection of $\tilde \gamma_2$ and $\si$.  Using the fact that the shortest distance between a point and a horocircle is always measured along a vertical geodesic, we easily see that $d_{\A}(\bar{B},\si)\leq d_{\A}(\bar{B},D)  \leq d_{\A}(\bar{B},B) \leq 1$, and thus the distance along the horocircle between $B$ and $D$ is at most $e^2$.  The points on the segment $\gamma_2(\omega_2(t))$ for $t \in [t_B,d]$ exactly lie on the horcyclic segment from $B$ to $D$ union the geodesic segment from $D$ to $\bar{B}$.  It is clear from the previous measurements that for any such point $q$ on this union, $d_{\A}(q,B) \leq 1+e^2$.

Lastly, consider $t \in [d,length(\tilde \xi_2)]$.  If $\tilde \gamma_2(\omega_2(t)) \in \Omega_p$, then the parametrization agrees with $\rho$ and it
immediately follows that
$d_{\A}(\xi_2(\omega_1(t)),\gamma_2(\omega_2(t))) \leq 1$. If
$\tilde \gamma_2(\omega_2(t)) \notin \Omega_p$, then the fact that
$d_{\A}(\xi_2(\omega_1(t)),\gamma_2(\omega_2(t))) \leq 2$
follows from lemma \ref{lemma:hyperbolic}.
\end{enumerate}

\noindent
Combining the above cases, we have shown that $d_{\A}(\xi_2(\omega_1(t)),\gamma_2(\omega_2(t))) \leq e^2 +1$. \end{proof}

We now use lemma \ref{lemma:bounded} to prove proposition \ref{prop:boundedwidth}.

\medskip
\noindent {\it Proof of proposition \ref{prop:boundedwidth}.} Let
$\alpha_1, \ \alpha_2 \in \Omega_p$ with $d(\alpha_1,\alpha_2)
\leq 1$. As above, the combing paths are given by
$\sigma_{\alpha_1} = \gamma_2 \circ \gamma_1$ and
$\sigma_{\alpha_2}=\xi_2 \circ \xi_1$.  Let $\psi \in R$ be the
reparametrization function described above satisfying
$d(\xi_1(\psi(t)),\gamma_1(t)) \leq 2 \log p$ for $t \in
[0,length(\gamma_1)]$.

Suppose that $\gamma_1=\xi_1$, so that the points $\alpha_1$ and $\alpha_2$ lie in a common hyperbolic plane.  Combining $\psi$ with the reparametrization functions $\omega_1$ and $\omega_2$ of $\xi_2$ and $\gamma_2$, to get $\bar{\omega_1}$ and $\bar{\omega_2}$ respectively, we see that $d(\sigma_{\alpha_1}(\bar{\omega_1}(t)),\sigma_{\alpha_2}(\bar{\omega_2}(t))) \leq e^2 + 2$ for all $t \in [0, length(\sigma_{\alpha_2})]$.

If $\gamma_1 \neq \xi_1$, so that the points $\alpha_1$ and
$\alpha_2$ lie in different hyperbolic planes within $\Omega_p$,
we recall that we use the product metric on $\HT$, which is
Lipschitz equivalent to the word metric on $\PZp$. The above
argument again bounds the asynchronous width between the paths by
the constant $\max \{ 2 \log p, \ e^2 + 2 \}$.

Thus we see that the asynchronous width of the combing of $\Omega_p$ defined above is bounded.
\qed

\begin{corollary}
\label{cor:Pcombing}
The group $\PZp$ has a combing with bounded asynchronous width.
\end{corollary}

\begin{proof}
From proposition \ref{prop:boundedwidth}, we know that the combing
constructed above of $\Omega_p$ has bounded asynchronous width.
Since $\Omega_p$ and $\PZp$ are quasi-isometric, choose a finite
generating set $X$ for $\PZp$ and let $f: \Omega_p \rightarrow
\PZp$ be a $(K,C)$-quasi-isometry, where $\PZp$ is viewed as
generated by $X$.   A different choice of generating set would
simply result in a different quasi-isometry.  Using the standard
connect-the-dots procedure (see, e.g., \cite{FS}), adapt $f$ so
that it is continuous.  As in theorem 3.6.4 of \cite{E+}, use $f$
to determine a combing of $\PZp$ in the generating set $X$ by
taking the images of the combing paths in $\Omega_p$. It follows
that for a constant $M$ depending on the quasi-isometry constants
$K$ and $C$ and the bound on the asynchronous width of the combing
above, the induced combing of $\PZp$ has asynchronous width
bounded by $M$.
\end{proof}

\section{The Dehn function of $\PZp$}
\label{sec:isoper}

We now prove the following theorem, using the combing constructed
in \S \ref{sec:combing}.

\begin{theorem}
\label{thm:mainthm} The Dehn function of $\PZp$ is exponential.
\end{theorem}

The proof of theorem \ref{thm:mainthm} consists of showing that
the Dehn function of $\PZp$ has both exponential upper and lower
bounds.

\subsection{An exponential lower bound}

We begin the proof that the Dehn function of $\PZp$ is exponential
by exhibiting an exponential lower bound, again using the space
$\Omega_p$. The Dehn function of $\Omega_p$ is well defined since
we can define Lipschitz loops in $\Omega_p$ and consider the area
of their fillings. It follows from \cite{BT} that the Dehn
functions of the group $\PZp$ and space $\Omega_p$ are equivalent,
in the sense described in \S \ref{sec:intro}. Thus if the Dehn
function of $\Omega_p$ has an exponential lower bound, then so
does the Dehn function of $\PZp$.

The space $\Omega_p$ has boundary components modelled on the
solvable Baumslag-Solitar group $BS(1,p^2)$, which we know has
exponential Dehn function.  A simple projection argument, given
below, easily shows that the Dehn function of $\Omega_p$, and
hence $\PZp$, has an exponential lower bound.

Let $X$ be the space $\HT$ with the interior of the horosphere
$\si$ based at $\infty$  removed. Since the Dehn function of
$BS(1,n)$, for integral $n > 1$, is exponential, there are
families of loops in $\si=X_n$, the $2$-complex defined in \S
\ref{sec:horospheres}, whose minimal filling is exponential in the
length of the loop. (See \cite{E+} for the precise definition of
these loops.)

\begin{lemma}
\label{lemma:lower} The Dehn function of $\PZp$ has an exponential
lower bound.
\end{lemma}

\begin{proof}
Let $f: S^1 \rightarrow \si$ be a Lipschitz loop in $\si$, with
$\gamma = f(S^1)$, whose minimal filling in $\si$ is exponential
in $n = length(\gamma)$. Let $\hat f: D^2 \rightarrow \Omega_p$ be
any filling of $\gamma$ lying in $\Omega_p$. In order to show that
the Dehn function of $\Omega_p$ is at least exponential, we
construct a projection $\pi: X \rightarrow \si$ which does not
increase the area of $F = \hat f(D^2)$ by more than a bounded
amount. For each $\A \subset \HT$, let the map $\pi$ project
upwards in $\Omega_p \cap \A$ along vertical geodesics to $\si
\cap \A$. In coordinates, if $h(t)$ denotes the height of the
point $t \in T_p$, then $\pi(x,y,t) = (x,p^{2h(t)} B,t)$ where $B$
was chosen in \S \ref{sec:omegap}.

Let $(dx,dy,dt)$ be the canonical coordinates for $\Omega_p$
viewed as a subspace of $\R^3$. If $S \subset \Omega_p$ is a
segment of the form $(x,y,[t_1,t_2])$, where $t_1$ and $t_2$ are
adjacent vertices in $T_p$ and $(x,y)$ is fixed, then $\pi(S)$ has
length $2 \log p$ in $\si$ whereas $S$ has length $1$ in
$\Omega_p$. \cite{T} Hence $\pi_*(dt) = Ldt$ where $L=2 \log p$.
So the area of $\pi(F)$, which is at least exponential in the
length of $\gamma$, is at most $L$ times the area of $F$. This
provides an exponential lower bound for the Dehn function of
$\Omega_p$. It follows from \cite{BT} that the Dehn functions of
$\Omega_p$ and $\PZp$ are equivalent, and thus the Dehn function
of $\PZp$ also satisfies an exponential lower bound.
\end{proof}

\subsection{An exponential upper bound}

To prove theorem \ref{thm:mainthm}, we must obtain an exponential
upper bound on the Dehn function of $\PZp$.  To do this, we apply
the following theorem of Bridson \cite{B} to the combing of $\PZp$
obtained in \S \ref{sec:combing}.

\begin{theorem}[\cite{B}, thms. 5.2(a), 6.1]
\label{thm:bridson} Suppose that $\sigma: G \rightarrow X^*$ is a
combing with asynchronous width $\phi(n)$ and length $L(n)$.  If
there exist constants $\alpha > 1$ and $C>0$ so that the
inequalities $L(n) \leq e^{Cn}$ and $\phi(n) < \alpha n$ for large
$n$, then there is a constant $k>0$ so that $e^{kn}$ is an
isoperimetric  function for some finite presentation of $G$.
\end{theorem}

\begin{lemma}
\label{lemma:upper} The Dehn function of $\PZp$ has an exponential
upper bound.
\end{lemma}

%\bigskip \noindent {\it Proof of theorem \ref{thm:mainthm}.}
%In \S \ref{sec:combing} we construct a combing
%of $\Omega_p$ in which the length of a combing path is at most exponential in the
%distance from the origin.  In addition, proposition \ref{prop:boundedwidth} proves that the
%combing has bounded asynchronous width.  Let $f: \Omega_p \rightarrow \PZp$ be a
%$(K,C)$-quasi-isometry, and use $f$ to determine a combing of the Caley graph of
%$\PZp$ by taking the images of the combing paths in $\Omega_p$.
%It follows that for slightly larger constants, which will depend on the
%constants $K$ and $C$
%of the quasi-isometry, this new combing also
%as boundedly asynchronous width and satisfies the same exponential bound on the
%engths of the combing paths.  Thus this new combing of $\PZp$ satisfies the confitions of
%theorem \ref{thm:bridson}, and we conclude
%that in some finite presentation, hence all finite presentations, of $\PZp$ have an
%exponential Dehn function.

\begin{proof}
In \S \ref{sec:combing}, we constructed a combing of $\Omega_p$
which satisfied $L(n) \leq e^{n}$ for all positive integers $n$.
If $f: \Omega_p \rightarrow \PZp$ is a $(K,C)$-quasi-isometry,
which has been adapted to be continuous, as in the proof of
corollary \ref{cor:Pcombing} we use $f$ to determine a combing of
$\PZp$.

Lemma \ref{lemma:length} states that the length function $L(n)$ of
the combing of $\Omega_p$ satisfies $L(n) \leq e^n$.  Thus, for a
constant $M = M(C,K)$, it follows that the length function $L(n)$
of the combing of $\PZp$ satisfies $L(n) \leq e^{Mn}$. From
corollary \ref{cor:Pcombing} we know that this combing of $\PZp$
has bounded asynchronous width. Thus both conditions of theorem
\ref{thm:bridson} are satisfied, and we conclude that in some
finite presentation, hence all finite presentations, $\PZp$ has an
exponential isoperimetric function, i.e. the Dehn function has an
exponential upper bound.
\end{proof}

%\subsection{The Dehn function of $\PZp$}
%\label{sec:proof}

Combining lemmas \ref{lemma:lower} and \ref{lemma:upper}, we
obtain the proof of theorem \ref{thm:mainthm}.  The following
corollary is immediate.

\begin{corollary}
$\PZp$ is not an automatic group.
\end{corollary}

\begin{proof}
It is shown in theorem 2.3.12 of \cite{E+} that an automatic group satisfies a quadratic
isoperimetric inequality, namely its Dehn function is bounded above by
a quadratic function.  It then follows from theorem \ref{thm:mainthm} that
$\PZp$ is not an automatic group.
\end{proof}

\end{document}